# APPROXIMATELY UNBIASED TESTS OF REGIONS USING MULTISTEP-MULTISCALE BOOTSTRAP RESAMPLING[1]

By Hidetoshi Shimodaira

*Tokyo Institute of Technology*

Approximately unbiased tests based on bootstrap probabilities are considered for the exponential family of distributions with unknown expectation parameter vector, where the null hypothesis is represented as an arbitrary-shaped region with smooth boundaries. This problem has been discussed previously in Efron and Tibshirani [*Ann. Statist.* **26** (1998) 1687–1718], and a corrected *p*-value with second-order asymptotic accuracy is calculated by the two-level bootstrap of Efron, Halloran and Holmes [*Proc. Natl. Acad. Sci. U.S.A.* **93** (1996) 13429–13434] based on the ABC bias correction of Efron [*J. Amer. Statist. Assoc.* **82** (1987) 171–185]. Our argument is an extension of their asymptotic theory, where the geometry, such as the signed distance and the curvature of the boundary, plays an important role. We give another calculation of the corrected *p*-value without finding the "nearest point" on the boundary to the observation, which is required in the two-level bootstrap and is an implementational burden in complicated problems. The key idea is to alter the sample size of the replicated dataset from that of the observed dataset. The frequency of the replicates falling in the region is counted for several sample sizes, and then the *p*-value is calculated by looking at the change in the frequencies along the changing sample sizes. This is the multiscale bootstrap of Shimodaira [*Systematic Biology* **51** (2002) 492–508], which is third-order accurate for the multivariate normal model. Here we introduce a newly devised multistep-multiscale bootstrap, calculating a third-order accurate *p*-value for the exponential family of distributions. In fact, our *p*-value is asymptotically equivalent to those obtained by the double bootstrap of Hall [*The Bootstrap and Edgeworth Expansion* (1992) Springer, New York] and the modified signed likelihood ratio of Barndorff-Nielsen [*Biometrika* **73** (1986) 307–322] ignoring $O(n^{-3/2})$ terms, yet the computation is less demanding and free from model specification.

Received November 2000; revised March 2004.

[1]Supported in part by Grant KAKENHI-14702061 from MEXT of Japan.

The algorithm is remarkably simple despite complexity of the theory behind it. The differences of the $p$-values are illustrated in simple examples, and the accuracies of the bootstrap methods are shown in a systematic way.

**1. Introduction.** We start with a simple example of Efron and Tibshirani (1998) to illustrate the issue to discuss. Let $X_1, \ldots, X_n$ be independent $p$-dimensional multivariate normal vectors with mean vector $\mu$ and covariance matrix identity $I_p$,

$$X_1, \ldots, X_n \sim N_p(\mu, I_p).$$

For given observed values $x_1, \ldots, x_n$, let us assume that we would like to know whether $\|\mu\|^2 = \mu_1^2 + \cdots + \mu_p^2 \leq 1$ or not. The problem is also described in a transformed variable $Y = \sqrt{n}\bar{X}$ with mean $\eta = \sqrt{n}\mu$, where $\bar{x} = (x_1 + \cdots + x_n)/n$ is the sample average. We have observed a $p$-dimensional multivariate normal vector $y$ having unknown mean vector $\eta$ and covariance matrix the identity,

$$(1.1) \qquad\qquad Y \sim N_p(\eta, I_p).$$

Then the null hypothesis we are going to test is $\eta \in \mathcal{R}$, with the spherical region

$$(1.2) \qquad\qquad \mathcal{R} = \{\eta : \|\eta\| \leq \sqrt{n}\,\}.$$

This problem is simple enough to give the exact answer. The frequentist confidence level, namely, the probability value ($p$-value) for the spherical null hypothesis is calculated as the probability of $\|Y\|^2$ being greater than or equal to the observed $\|y\|^2$ assuming that $\eta$ is on the boundary $\partial \mathcal{R} = \{\eta : \|\eta\| = \sqrt{n}\,\}$ of $\mathcal{R}$. The exact $p$-value is easily calculated knowing that $\|Y\|^2$ is distributed as the chi-square distribution with degrees of freedom $p$ and noncentrality $\|\eta\|^2$.

In this paper we are going to remove two restrictions in the above problem for generalization. (i) The underlying probability model for $Y$ is the exponential family of distributions, instead of the multivariate normal model; we denote the density function with the expectation parameter $\eta$ as

$$(1.3) \qquad\qquad Y \sim f(y; \eta).$$

(ii) The null hypothesis will be represented as an arbitrarily-shaped region $\mathcal{R}$ with smooth boundaries, instead of the spherical region. The surface of $\partial \mathcal{R}$ may be represented as the Taylor series with coefficients $d^{ab}, e^{abc}, \ldots$

$$(1.4) \qquad \Delta\eta_p = -d^{ab}\Delta\eta_a \Delta\eta_b - e^{abc}\Delta\eta_a \Delta\eta_b \Delta\eta_c + \cdots$$

in the local coordinates $(\Delta\eta_1, \ldots, \Delta\eta_p)$ by taking the origin at a point on $\partial \mathcal{R}$ and rotating the axes properly. The summation convention such as $d^{ab}\Delta\eta_a\Delta\eta_b = \sum_{a=1}^{p-1}\sum_{b=1}^{p-1} d^{ab}\Delta\eta_a\Delta\eta_b$ will be used, where the indices $a, b, \ldots$ may run through $1, \ldots, p-1$ and $i, j, \ldots$ may run though $1, \ldots, p$ when used as subscripts or superscripts for $p$-dimensional vectors. The axes are taken so that $\Delta\eta_1, \ldots, \Delta\eta_{p-1}$ are for the tangent space of the surface, and $\Delta\eta_p$ is for its orthogonal space taken positive in the direction pointing away from $\mathcal{R}$. This general setting is the "problem of regions" discussed previously in Efron and Tibshirani ([1998]), and our argument is an extension of their asymptotic theory, where the geometry, such as the signed distance and the curvature of the boundary, plays an important role.

Since the exact $p$-value is available only for special cases, we will discuss several bootstrap methods to calculate approximate $p$-values from $y$ under the assumptions (i) and (ii) above. Let $\alpha$ denote a specified significance level, and $\hat{\alpha}(y)$ denote an approximate $p$-value. A large value of $\hat{\alpha}(y)$ may indicate evidence to support the null hypothesis $\eta \in \mathcal{R}$. On the other hand, if $\hat{\alpha}(y) < \alpha$ is observed, then we reject the null hypothesis and conclude that $\eta \notin \mathcal{R}$. The hypothesis test of $\mathcal{R}$ is said to be *unbiased* if the rejection probability is equal to $\alpha$ whenever $\eta \in \partial\mathcal{R}$. The approximate $p$-value is said to be $k$th order accurate if the asymptotic bias is of order $O(n^{-k/2})$, that is,

$$(1.5) \qquad \Pr\{\hat{\alpha}(Y) < \alpha; \eta\} = \alpha + O(n^{-k/2}), \qquad \eta \in \partial\mathcal{R},$$

holds for $0 < \alpha < 1$. For sufficiently large $n$, approximately unbiased $p$-values of higher-order accuracy are considered to be better than those of lower-order accuracy.

We will not specify the probabilistic model or the shape of the region explicitly in the calculation of the $p$-value, but only assume that a mechanism is available to us for generating the bootstrap replicates and identifying whether the outcomes are in the region or not. This setting is important for complicated practical applications, where the exact $p$-value is not available and, thus, bootstrap methods are used for approximation. The phylogenetic tree selection discussed in Efron, Halloran and Holmes ([1996]) and Shimodaira ([2002]) is a typical case; the history of evolution represented as a tree is inferred by a model-based clustering of the DNA sequences of organisms, where we are given complex computer software for inferring the tree from a dataset. For calculating $p$-values of the hypothetical evolutionary trees, we can easily run bootstrap simulations, although computationally demanding, by repeatedly applying the software to replicated datasets.

We confine our attention to the parametric bootstrap of continuous random vectors for mathematical simplicity. We also assume that the boundary of the region is a smooth surface. In practical applications, however, it is often the case that the nonparametric bootstrap is employed, the random vector is discrete and the boundary is nonsmooth. Regions with nonsmooth boundaries, in particular, may lead to serious difficulty as discussed in Perlman and Wu ([1999](#), [2003](#)). Further study is needed to bridge these gaps between the theory and practice.

The frequency of the bootstrap replicates falling in the region, namely, the bootstrap probability, has been used widely since its application to phylogenetic tree selection in Felsenstein ([1985](#)). This is also named "empirical strength probability" of $\mathcal{R}$ in Liu and Singh ([1997](#)), where a modification for nonsmooth boundary is discussed as well. The bootstrap probability is, however, biased as an approximation to the exact $p$-value and, thus, the *two-level bootstrap* of Efron, Halloran and Holmes ([1996](#)) and Efron and Tibshirani ([1998](#)) is developed to improve the accuracy. Under the assumptions (i) and (ii) above, the two-level bootstrap calculates a second-order accurate $p$-value, whereas the bootstrap probability is only first-order accurate.

The bias of the bootstrap probability mainly arises from the curvature of $\partial \mathcal{R}$. The two-level bootstrap estimates the curvature for bias correction, where the curvature is estimated by generating second-level replicates around $\hat{\eta}(y)$. Here $\hat{\eta}(y)$ denotes the maximum likelihood estimate for $\eta$ restricted to $\partial \mathcal{R}$. $\hat{\eta}(y)$ is the nearest point on $\partial \mathcal{R}$ to $y$ for ([1.1](#)). For the spherical region, $\hat{\eta}(y) = \sqrt{n}y/\|y\|$ is easily obtained, but $\hat{\eta}(y)$ must be obtained by numerical search in general, leading to an implementational burden in complex problems. This motivated our development of a new method.

The *multiscale bootstrap* is developed in Shimodaira ([2002](#)) to calculate another bias corrected $p$-value. It does not require $\hat{\eta}(y)$. Instead, the bootstrap probabilities are calculated for sets of bootstrap replicates with several sample sizes which may differ from that of the observed data. This, in effect, alters the scale parameter of the replicates (Figure [1](#)). The key idea is to estimate the curvature from the change in the bootstrap probabilities along varying sample sizes. The corrected $p$-value is third-order accurate for any arbitrarily-shaped region with smooth boundaries under the multivariate normal model. The normality assumption is not as restrictive as it might look at first, because the procedure is transformation-invariant and should work fine if there exists a transformation from the dataset to the normal $Y$ and if the null hypothesis is represented as a region of $\eta$. We do not have to know what the transformation is. However, it becomes only first-order accurate if there is no such transformation to ([1.1](#)) but only one to ([1.3](#)).

The multiscale bootstrap can be used easily for complex problems. It is as easy as the usual bootstrap. We only have to change the sample size of the bootstrap replicates, and apply a regression fit to the bootstrap probabilities. The bias corrected $p$-value is calculated from the slope of the regression curve (Figure [2](#)). This procedure is implemented in computer software [Shimodaira and Hasegawa ([2001](#))] for phylogenetic tree selection, and is also applied to gene network estimation from microarray expression profiles [Kamimura et al. ([2003](#))]. In these applications, the multiscale bootstrap can calculate the $p$-values for many related hypotheses at the same time; we do not have to run time-consuming bootstrap simulations separately for these hypotheses. For example, biologists are interested in the monophyletic hypothesis that some specified species constitute a cluster in the phylogenetic tree, and there are many such hypotheses for groups of species. The bootstrap probabilities for these hypotheses are obtained at the same time from a single run of bootstrap simulation for each scale. We only have to apply the regression fit separately to the multiscale bootstrap probabilities of each hypothesis.

In this paper we provide the theoretical foundation of the multiscale bootstrap, and introduce a newly devised *multistep-multiscale bootstrap* resampling. This method calculates an approximately unbiased $p$-value with third-order asymptotic accuracy under the assumptions (i) and (ii). The previously developed method of Shimodaira (2002) corresponds to a special case of the new method, that is, the one-step multiscale bootstrap.

For explaining the bootstrap methods, a rather intuitive argument is given in Sections 2 to 6 using simple examples. A more formal argument is given in Section 7, and the technical details are given in a supporting document [Shimodaira (2004)]. We introduce a *modified signed distance*, and give a unified approach to the asymptotic analysis of the bootstrap methods using Edgeworth series, as well as the tube formula of Weyl (1939). Third-order accuracy is also shown there for the $p$-value computed by the modified signed likelihood ratio [Barndorff-Nielsen (1986)], which requires the analytic expression of the likelihood function, and for the $p$-value computed by the double bootstrap [Hall ([1992](#))], which requires a huge number of replicates, as well as computation of $\hat{\eta}(y)$. The multistep-multiscale bootstrap method requires only the bootstrap mechanism for generating replicates around $y$, inheriting the simplicity from the one-step multiscale bootstrap. The price for higher-order accuracy and simpler implementation is a large number of replicates, which can be as large as that of the double bootstrap. These three $p$-values are, in fact, shown to be equivalent ignoring $O(n^{-3/2})$ terms.

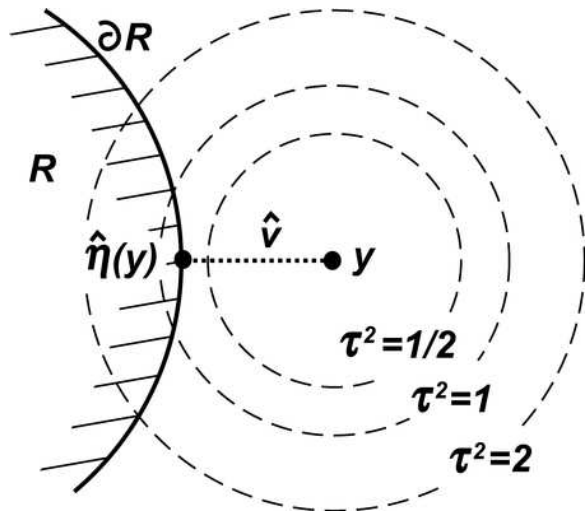

Fig. 1.  *Multiscale bootstrap. The three circles with dashed lines indicate the conditional distributions of the bootstrap replicates with mean $y$ and scales $\tau = 1/\sqrt{2}, 1, \sqrt{2}$. In this particular configuration, the bootstrap probability may increase by halving the sample size to alter $\tau = 1$ to $\sqrt{2}$, and may decrease by doubling the sample size to alter $\tau = 1$ to $1/\sqrt{2}$.*

Our argument may not be justified unless the assumptions (i) and (ii) hold. We are not sure yet how robust the multistep-multiscale bootstrap method is under misspecifications of the exponential family model. It is shown at the end of Section [4](#), however, that the one-step method adjusts the bias halfway, though not completely, under misspecifications of the normal model. A simulation study in Shimodaira ([2002](#)) shows that the bias of the one-step method under the normal model is very small even if the boundary is piecewise smooth, but the bias becomes larger as $\eta$ moves closer to nonsmooth points on the boundary.

**2. Two-level bootstrap resampling.** Although our ultimate goal is to get rid of the normal assumption, we use normality in this section to illustrate the bootstrap methods, and besides ([1.1](#)), we also assume ([1.2](#)). For given observed value $\bar{x}$, we consider the parametric bootstrap resampling

$$X_1^*, \dots, X_{n_1}^* \sim N_p(\bar{x}, I_p).$$

Typically, the sample size $n_1$ of the replicated dataset should be equal to $n$, but we reserve the generality of using any value for $n_1$. The scaling factor of the bootstrap, $\tau_1 = \sqrt{n/n_1}$, will be altered later in the multiscale bootstrap. Once we specify $\tau_1$, we may generate $B$, say 10,000, replicated datasets, and compute the average $\overline{X}^* = (X_1^* + \dots + X_{n_1}^*)/n_1$ for each replicate. A large value of the frequency that $\|\overline{X}^*\|^2 \leq 1$ holds in the replicates may indicate a high chance of the null hypothesis $\|\mu\|^2 \leq 1$ being correct. This is also described in a transformed variable $Y^* = \sqrt{n}\overline{X}^*$. For given observed value $y$, we consider the parametric bootstrap resampling

$$(2.1) \qquad\qquad Y^* \sim N_p(y, \tau_1^2 I_p),$$

and the bootstrap probability with scale $\tau_1$ is denoted by

$$\tilde{\alpha}_1(y, \tau_1) = \Pr\{Y^* \in \mathcal{R}; y, \tau_1\},$$

where the index 1 indicates the "one-step" bootstrap in connection with $\tilde{\alpha}_2$ and $\tilde{\alpha}_3$ defined later, as shown in Table [1](#). $\tilde{\alpha}_1$ is estimated by the frequency of $Y^* \in \mathcal{R}$ from the $B$ bootstrap replicates with the binomial variance $\tilde{\alpha}_1(1 - \tilde{\alpha}_1)/B$.

Let us consider a numerical example with

$$(2.2) \qquad p = 4, \qquad n = 10, \qquad \|\bar{x}\|^2 = 2.680.$$

Although $\|\bar{x}\|^2 > 1$, we are not sure if $\|\mu\|^2 \leq 1$ holds or not. The frequentist confidence level for the null hypothesis is given by the exact $p$-value, which we will denote by $\hat{\alpha}_\infty(y)$, or simply $\hat{\alpha}_\infty$ for brevity sake. In this numerical example, the value of $\|\bar{x}\|^2$ is, in fact, chosen to make $\hat{\alpha}_\infty(y) = 0.05$. $\hat{\alpha}_\infty$ may be approximated by the bootstrap probability with $\tau_1 = 1$, denoted by

$$\hat{\alpha}_0(y) = \tilde{\alpha}_1(y, 1).$$

This turns out to be $\hat{\alpha}_0(y) = 0.0085$, showing $\hat{\alpha}_0$ is not a very good approximation to $\hat{\alpha}_\infty$. Here the problem is so simple that $\hat{\alpha}_0(y)$, as well as $\hat{\alpha}_\infty(y)$, can be computed numerically from the noncentral chi-square distribution function. If the bootstrap resampling with $B = 10{,}000$, say, is used for $\hat{\alpha}_0$, the standard error becomes $0.0009$.

A modification of $\hat{\alpha}_0$ is developed based on the geometric theory in Efron, Halloran and Holmes ([1996](#)) and Efron and Tibshirani ([1998](#)) to improve the accuracy of the approximation to $\hat{\alpha}_\infty$. The idea is to compute $\hat{\alpha}_0(\hat{\eta}(y))$ by generating the second-level replicates around $\hat{\eta}(y)$ for estimating the curvature of the surface $\partial\mathcal{R}$. When the surface of $\partial\mathcal{R}$ is flat, $\hat{\alpha}_0(\hat{\eta}(y)) = \frac{1}{2}$. It becomes smaller/larger than $\frac{1}{2}$ when the surface is curved toward/away from $\mathcal{R}$. Let $z$ denote a generic symbol for the $z$-value corresponding to a $p$-value $\alpha$ with relation $z = -\Phi^{-1}(\alpha)$, where $\Phi^{-1}(\cdot)$ is the inverse of the standard normal distribution function $\Phi(\cdot)$. For example, we may write $\hat{z}_0(y) = -\Phi^{-1}(\hat{\alpha}_0(y))$. The *ABC conversion formula* of Efron ([1987](#)) and DiCiccio and Efron ([1992](#)) is

$$(2.3) \qquad \hat{z}_{\mathrm{abc}}(y) = \frac{\hat{z}_0(y) - \hat{z}_0(\hat{\eta}(y))}{1 - \hat{a}(\hat{z}_0(y) - \hat{z}_0(\hat{\eta}(y)))} - \hat{z}_0(\hat{\eta}(y)),$$

TABLE 1
*Bootstrap probabilities and corrected p-values*

| Symbol | Section | Description |
|---|---|---|
| $\tilde{\alpha}_1(y, \tau_1)$ | 2 | Bootstrap probability |
| $\hat{\alpha}_\infty(y)$ | 2 | Exact $p$-value[*] |
| $\hat{\alpha}_0(y)$ | 2 | Bootstrap probability ($\tau_1 = 1$) |
| $\hat{\alpha}_{\mathrm{abc}}(y)$ | 2 | Two-level bootstrap corrected $p$-value |
| $\hat{\alpha}_1(y)$ | 3 | Multiscale bootstrap corrected $p$-value |
| $\tilde{\alpha}_2(y, \tau_1, \tau_2)$ | 4 | Two-step bootstrap probability |
| $\hat{\alpha}_2(y)$ | 4 | Two-step multiscale bootstrap corrected $p$-value |
| $\tilde{\alpha}_3(y, \tau_1, \tau_2, \tau_3)$ | 5 | Three-step bootstrap probability |
| $\hat{\alpha}_3(y)$ | 5 | Three-step multiscale bootstrap corrected $p$-value |

[*]A third-order accurate $p$-value in Section 7.

where $\hat{z}_{abc}(y)$, $\hat{z}_0(y)$, and $\hat{z}_0(\hat{\eta}(y))$ are denoted $\widehat{Z}$, $\widetilde{Z}$, and $\hat{z}_0$, respectively, in the notation of equation (6.6) of Efron and Tibshirani (1998). The corrected $p$-value for the two-level bootstrap is then defined by $\hat{\alpha}_{abc}(y) = \Phi(-\hat{z}_{abc}(y))$. The *acceleration constant* $\hat{a}$, characterizing the probabilistic model, is known to be $\hat{a} = 0$ for the normal model. $\hat{a}$ may also be estimated using the second-level bootstrap for (1.3); for details we refer to Efron, Halloran and Holmes (1996). Note that the sign in front of $\hat{a}$ in (2.3) is reversed from that of equation (6.6) of Efron and Tibshirani (1998), because the $\Delta\eta_p$-axis is taking the opposite direction here.

The $p$-values for the numerical example of (2.2) are

$$\hat{\alpha}_0(y) = 0.0085, \qquad \hat{\alpha}_0(\hat{\eta}(y)) = 0.315,$$

$$\hat{\alpha}_{abc}(y) = 0.0775, \qquad \hat{\alpha}_\infty(y) = 0.05.$$

We observe that $\hat{\alpha}_{abc}$ shows great improvement over $\hat{\alpha}_0$ to approximate $\hat{\alpha}_\infty$. This improvement is also confirmed in the asymptotic argument. It has been shown in Efron and Tibshirani (1998) that $k = 1$ for $\hat{\alpha}_0$, and $k = 2$ for $\hat{\alpha}_{abc}$ under (1.3) and (1.4).

**3. Multiscale bootstrap resampling.** Here we continue to use the normal model (1.1) for the argument of the corrected $p$-value in this section. The bootstrap probability changes if the replicate sample size changes. When we alter $n_1 = 10$ to $n_1 = 3$ for the numerical example of (2.2), or equivalently alter the scale $\tau_1 = 1$ to $\tau_1 = \sqrt{10/3}$, we observe that $\hat{\alpha}_1(y, 1) = 0.0085$ changes to $\hat{\alpha}_1(y, \sqrt{10/3}) = 0.0359$. In the multiscale bootstrap, $\hat{\alpha}_1(y, \tau_1)$ is computed for several values of $\tau_1 = \sqrt{n/n_1}$. For example, instead of $n = 10$, we use the following five $n_1$ values:

$$(3.1) \qquad\qquad n_1 = 3, 6, 10, 15, 21,$$

and compute the corresponding bootstrap probabilities

$$(3.2) \qquad \tilde{\alpha}_1(y, \tau_1) = 0.0359, 0.0205, 0.0085, 0.0028, 0.0008.$$

These values, as well as those for other parameter settings, are shown in Figure 2 by plotting the $z$-value along the inverse of the scale. The horizontal axis is $1/\tau_1 = \sqrt{n_1/n} = 0.55, 0.78, 1, 1.23, 1.45$, and the vertical axis is $\tilde{z}_1(y, \tau_1) = -\Phi^{-1}(\tilde{\alpha}_1(y, \tau_1)) = 1.80, 2.04, 2.39, 2.77, 3.17$.

Figure 2 shows these values along with a regression fit. This is obtained by fitting a regression model with explanatory variables $1/\tau_1$ and $\tau_1$,

$$(3.3) \qquad\qquad \tilde{z}_1(y, \tau_1) \approx \hat{v}/\tau_1 + \hat{c}\tau_1,$$

to the plot, where $\hat{v}$ and $\hat{c}$ are the regression coefficients estimated as

$$(3.4) \qquad\qquad \hat{v} = 2.002, \qquad \hat{c} = 0.385$$

for the plot of (3.2). We observe that the regression fit agrees with the plots very well for the cases in Figure 2. The regression model (3.3) has been justified in Shimodaira ([2002]) under ([1.1]) and ([1.4]); we will use "$\approx$" to indicate that equality holds up to $O(n^{-1})$ terms with the error of order $O(n^{-3/2})$. The regression model with explanatory variables $1/\tau_1$ and $\tau_1$ will be justified later, in fact, under ([1.3]) and ([1.4]) as seen in ([7.15]), although the following interpretation of the coefficients should be modified accordingly.

A simple geometric interpretation can be given to the regression coefficients under ([1.1]) and ([1.4]). Efron and Tibshirani ([1998]) have shown a formula equivalent to

$$\hat{z}_0(y) \approx \hat{v} + \hat{c}, \tag{3.5}$$

where $\hat{v}$ and $\hat{c}$ correspond to $x_0$ and $\hat{d}_1 - x_0\hat{d}_2$, respectively, in their equation ([2.19]). $\hat{v}$ is the *signed distance* of Efron ([1985]), defined as the distance from $y$ to $\partial\mathcal{R}$ with a positive/negative sign when $y$ is outside/inside of $\mathcal{R}$. Thus, $\hat{v} = \pm\|y - \hat{\eta}(y)\|$ measures evidence of the null hypothesis being wrong. $\hat{c}$ is related to the $(p-1) \times (p-1)$ matrix $\hat{d}^{ab}$ measuring the curvature of $\partial\mathcal{R}$ at $\hat{\eta}(y)$; $\hat{d}^{ab}$ is defined as $d^{ab}$ in ([1.4]) by making the local coordinates orthonormal at $\hat{\eta}(y)$. In our notation, $\hat{c} = \hat{d}_1 - \hat{v}\hat{d}_2$, where $\hat{d}_1 = \hat{d}^{aa}$ is the trace of $\hat{d}^{ab}$, and $\hat{d}_2 = (\hat{d}^{ab})^2 = \sum_{a=1}^{p-1}\sum_{b=1}^{p-1}(\hat{d}^{ab})^2$ is that for the squared matrix. When $\partial\mathcal{R}$ is flat at $\hat{\eta}(y)$, $\hat{d}^{ab} = 0$ and, thus, $\hat{c} = 0$. $\hat{v}$, $\hat{d}_1$ and $\hat{d}_2$ are transformation-invariant functions of $y$ calculated from the shape of the boundary and the density function of $Y$; they are referred to as geometric

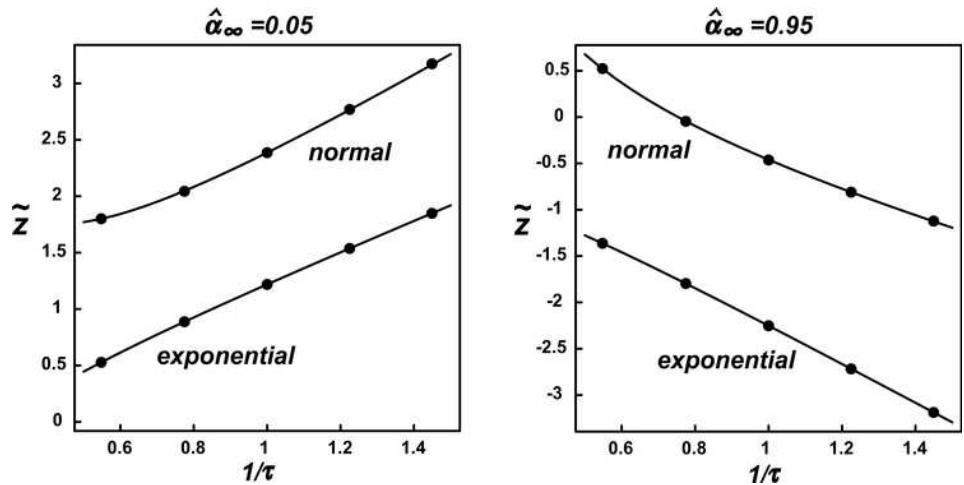

FIG. 2. *Plots of the z-value of the multiscale bootstrap probability along the inverse of the scale $\tau$ for the normal example ($p = 4$) of Section [2] and the exponential example ($p = 1$) of Section [4]. Parameter values are chosen so that the exact p-value is either 0.05 (left panel) or 0.95 (right panel). The curves are drawn by the regression model of equation ([3.3]).*

quantities here. Under (1.1) and (1.2) these quantities are

$$(3.6) \qquad \hat{v} = \|y\| - \sqrt{n}, \qquad \hat{d}_1 = \frac{p-1}{2\sqrt{n}}, \qquad \hat{d}_2 = \frac{p-1}{4n}.$$

This computes directly,

$$(3.7) \qquad \hat{v} = 2.015, \qquad \hat{c} = 0.323$$

for (2.2), showing good agreement with those computed indirectly from the multiscale bootstrap. $\hat{v}$ and $\hat{c}$ in (3.4) are actually estimating those in (3.7), thus, it would be appropriate to denote the former as $\hat{\hat{v}}$ and $\hat{\hat{c}}$, although we do not make the notational distinction. This estimation is third-order accurate, since the regression model (3.3) holds for (3.7) with error of $O(n^{-3/2})$.

Considering that $\hat{v}$ and $\hat{c}$ are functions of $y$, we may define a statistic

$$(3.8) \qquad \hat{z}_1(y) = \hat{v} - \hat{c}.$$

This is equivalent to the pivot statistic of Efron (1985), and $\Pr\{\hat{z}_1(Y) \leq x; \eta\} \approx \Phi(x)$ for $\eta \in \partial\mathcal{R}$ under (1.1) and (1.4); see equation (2.16) of Efron and Tibshirani (1998). Thus, a third-order accurate $p$-value is defined by $\hat{\alpha}_1(y) = \Phi(-\hat{z}_1(y))$. We can compute $\hat{\alpha}_1(y)$ using $\hat{v}$ and $\hat{c}$ obtained from the multiscale bootstrap. For the example of (2.2),

$$\hat{\alpha}_1(y) = \Phi(-2.002 + 0.385) = 0.0529,$$

showing an improvement over $\hat{\alpha}_{abc}(y) = 0.0775$ to approximate $\hat{\alpha}_\infty(y) = 0.05$. The index of $\hat{\alpha}_1$ indicates the "one-step" bootstrap as similarly for $\tilde{\alpha}_1$.

It is interesting to note that we can also read off the values of $\hat{z}_1(y)$ from Figure 2. The differentiation of (3.3) with respect to $1/\tau_1$ is

$$\frac{\partial \tilde{z}_1(y, \tau_1)}{\partial(1/\tau_1)} \approx \hat{v} - \hat{c}\tau_1^2,$$

and the slope of the regression curve at $1/\tau_1 = 1$ gives $\hat{z}_1(y)$. The corrected $p$-value $\hat{\alpha}_1$ is essentially obtained from the change of the bootstrap probability in the multiscale bootstrap.

**4. Two-step multiscale bootstrap resampling.** The one-step multiscale bootstrap described in Section 3 calculates a very accurate $p$-value for the arbitrarily-shaped region if there exists a transformation from the dataset to the normal model. However, it can be inaccurate if such a transformation does not exist even approximately. This restriction essentially comes from the fact that the covariance matrix of $y$ in (1.1) is constant with respect to $\eta$. The acceleration constant $\hat{a}$ of the ABC formula measures the rate of change in the covariance matrix, and $\hat{a}$ is assumed zero in the derivation of (3.8). Here we introduce the *two-step multiscale bootstrap* for estimating $\hat{a}$ to improve the accuracy of the one-step multiscale bootstrap.

A breakdown of the one-step multiscale bootstrap method is illustrated in the following example. Let $X_1, \ldots, X_n$ be one-dimensional independent exponential random variables with mean $\mu$,

$$X_1, \ldots, X_n \sim \exp(-x/\mu - \log \mu),$$

and let the null hypothesis of interest be $\mu \leq 1$. The exact $p$-value is calculated by knowing that a transformed variable $Y = \sqrt{nX}$ is distributed as Gamma with shape $n$ and mean $\eta = \sqrt{n}\mu$. We consider a numerical example with

$$(4.1) \qquad p = 1, \qquad n = 10, \qquad \bar{x} = 1.571,$$

so that $\hat{\alpha}_\infty(y) = 0.05$. The multiscale bootstrap probabilities for the five $n_1$ values in (3.1) are computed as

$$(4.2) \qquad \tilde{\alpha}_1(y, \tau_1) = 0.2990, 0.1875, 0.1115, 0.0622, 0.0322,$$

and the regression coefficients of (3.3) are estimated as $\hat{v} = 1.328, \hat{c} = -0.110$. Then the corrected $p$-value is computed as

$$(4.3) \qquad \hat{\alpha}_1(y) = \Phi(-1.328 - 0.110) = 0.0753.$$

Although this is an improvement over $\hat{\alpha}_0(y) = 0.112$, it is not as good as in the normal example above. The pivot (3.8) is not justified under (1.3) in general, and $\hat{\alpha}_1(y)$ is, in fact, only first-order accurate for the exponential example.

The two-step multiscale bootstrap is employed simply to generate a second-step replicate from every first-step replicate. Let us denote the conditional density of the first-step bootstrap replicate $Y^* = \sqrt{nX^*}$ as

$$(4.4) \qquad Y^* \sim f(y^*; y, \tau_1),$$

given mean $y = \sqrt{nX}$ and scale $\tau_1$ under (1.3), which reduces to $f(y^*; y, 1) = f(y^*; y)$ when $\tau_1 = \sqrt{n/n_1}$ is unity. This becomes (2.1) for (1.1), and Gamma with shape $n_1$ and mean $y$ for the exponential example. We generate a second-step replicate $Y^{**}$ for each $y^*$. The conditional density of $Y^{**}$ given $y^*$ takes the same form as (4.4), but with scale parameter $\tau_2 = \sqrt{n/n_2}$;

$$(4.5) \qquad Y^{**} \sim f(y^{**}; y^*, \tau_2).$$

For the normal example, (4.5) is equivalent to generating

$$X_1^{**}, \ldots, X_{n_2}^{**} \sim N_p(\bar{x}^*, I_p)$$

for given $\bar{x}^*$, and using the transformed variable $Y^{**} = \sqrt{nX^{**}}$. The two-step bootstrap probability with a pair of scales $(\tau_1, \tau_2)$ is then defined by

$$\tilde{\alpha}_2(y, \tau_1, \tau_2) = \Pr\{Y^{**} \in \mathcal{R}; y, \tau_1, \tau_2\}$$

$$= \int \tilde{\alpha}_1(y^*, \tau_2) f(y^*; y, \tau_1) \, dy^*,$$

where the integration is taken over the range of the components. We can write $\tilde{\alpha}_1(y, \tau_1) = \tilde{\alpha}_2(y, \tau_1, 0)$, because the conditional density of $Y^{**}$ converges to the point mass at $y^*$ by taking the limit $\tau_2 \to 0$. The two-step bootstrap might look similar to the double bootstrap of Hall (1992), but they are very different. We should generate thousands of $Y^{**}$ for given $y^*$ in the double bootstrap, but only one $Y^*$ in the two-step bootstrap.

Let us consider two $n_2$ values,

$$(4.6) \qquad n_2 = 6, 15,$$

for the normal example with parameter values (2.2). The two-step bootstrap probabilities are, for example,

$$\tilde{\alpha}_2(y, \sqrt{\tfrac{10}{6}}, \sqrt{\tfrac{10}{6}}) = 0.0359, \qquad \tilde{\alpha}_2(y, \sqrt{\tfrac{10}{10}}, \sqrt{\tfrac{10}{15}}) = 0.0205.$$

Of course, they give $\tilde{\alpha}_1(y, \sqrt{\tfrac{10}{3}})$ and $\tilde{\alpha}_1(y, \sqrt{\tfrac{10}{6}})$, respectively, in (3.2), because

$$\tilde{\alpha}_2(y, \tau_1, \tau_2) = \tilde{\alpha}_1(y, \sqrt{\tau_1^2 + \tau_2^2})$$

for (1.1). For the exponential example with parameter values (4.1), however,

$$\tilde{\alpha}_2(y, \sqrt{\tfrac{10}{6}}, \sqrt{\tfrac{10}{6}}) = 0.3063, \qquad \tilde{\alpha}_2(y, \sqrt{\tfrac{10}{10}}, \sqrt{\tfrac{10}{15}}) = 0.1866$$

are different, though very slightly, from $\tilde{\alpha}_1(y, \sqrt{\tfrac{10}{3}}) = 0.2990$ and $\tilde{\alpha}_1(y, \sqrt{\tfrac{10}{6}}) = 0.1875$, respectively, in (4.2). The difference of $\tilde{\alpha}_2(y, \tau_1, \tau_2)$ from $\tilde{\alpha}_1(y, \sqrt{\tau_1^2 + \tau_2^2})$ for (1.3) is explained by

$$(4.7) \qquad \tilde{z}_2(y, \tau_1, \tau_2) - \tilde{z}_1(y, \sqrt{\tau_1^2 + \tau_2^2}) \doteq \frac{\hat{a}\tau_1^2 \tau_2^2 (\hat{v}^2 - (\tau_1^2 + \tau_2^2))}{(\tau_1^2 + \tau_2^2)^{5/2}}.$$

We will use "$\doteq$" to indicate that equality holds up to $O(n^{-1/2})$ terms with error of order $O(n^{-1})$. Formula (4.7) and a revised regression model

$$(4.8) \qquad \tilde{z}_1(y, \tau_1) \doteq \frac{\hat{v} - 2\hat{a}\hat{v}^2}{\tau_1} + (\hat{d}_1 - \hat{a})\tau_1$$

for (1.3) are consequences of a more general argument with third-order accuracy shown in Section 7.

The key idea in the two-step multiscale bootstrap is to estimate $\hat{a}$ by looking at the difference of $\tilde{\alpha}_2(y, \tau_1, \tau_2)$ from $\tilde{\alpha}_1(y, \sqrt{\tau_1^2 + \tau_2^2})$. Once we compute $\tilde{\alpha}_1(y, \tau_1)$ and $\tilde{\alpha}_2(y, \tau_1, \tau_2)$ for several values of $(\tau_1, \tau_2)$ by the one-step and two-step multiscale bootstrap, we can estimate $\hat{v}$, $\hat{d}_1$ and $\hat{a}$ by fitting (4.7) and (4.8) to the observed bootstrap probabilities. A second-order accurate $p$-value, denoted $\hat{\alpha}_2(y)$, is then computed by using the estimated geometric quantities in the $z$-value

$$(4.9) \qquad \hat{z}_2(y) \doteq \hat{v} - \hat{d}_1 + \hat{a}(1 - \hat{v}^2).$$

This expression is shown to be equivalent to (2.3) up to $O(n^{-1/2})$ terms by using (4.8); $\hat{z}_0(y) \doteq \hat{v} + \hat{d}_1 - \hat{a}(1 + 2\hat{v}^2)$ and $\hat{z}_0(\hat{\eta}(y)) \doteq \hat{d}_1 - \hat{a}$. In the next section we will describe a procedure based on the above idea, as well as its refined version with third-order accuracy.

It follows from (4.8) that the one-step multiscale bootstrap estimates $\hat{v} - 2\hat{a}\hat{v}^2$ and $\hat{d}_1 - \hat{a}$ for the coefficients $\hat{v}$ and $\hat{c}$, respectively, under (1.3). Thus, $\hat{z}_1(y) \doteq \hat{v} - \hat{d}_1 + \hat{a}(1 - 2\hat{v}^2) \doteq \hat{z}_2(y) - \hat{a}\hat{v}^2$, as well as $\hat{z}_0(y) \doteq \hat{z}_2(y) + 2\hat{d}_1 - 2\hat{a} - \hat{a}\hat{v}^2$, is first-order accurate in general. Since the difference $\hat{z}_2(y) - \hat{z}_1(y) \doteq \hat{a}\hat{v}^2$ does not involve $\hat{d}_1$, the one-step method adjusts the bias resulting from the curvature even if the normal model is misspecified.

**5. Three-step multiscale bootstrap resampling.** We may repeat "stepping" to obtain multistep-multiscale bootstrap probabilities so that we might be able to compute higher-order accurate $p$-values. This is the case, in fact, for going one step further, although the results are not known for yet further stepping. We introduce the *three-step multiscale bootstrap* for computing a third-order accurate $p$-value, denoted $\hat{\alpha}_3(y)$, under (1.3) and (1.4). In the following argument, we first describe the procedure to compute $\hat{\alpha}_2(y)$, which helps understand that for $\hat{\alpha}_3(y)$.

The expression for $\hat{z}_2(y, \tau_1, \tau_2)$ is obtained from (4.7) by substituting $\sqrt{\tau_1^2 + \tau_2^2}$ for $\tau_1$ in (4.8). This is also expressed as

$$(5.1) \qquad \tilde{z}_2(y, \tau_1, \tau_2) \doteq \zeta_2(\hat{\gamma}_1, \hat{\gamma}_2, \hat{\gamma}_3, \tau_1, \tau_2),$$

where the function $\zeta_2$ on the right-hand side is defined by

$$(5.2) \qquad \zeta_2(\gamma_1, \gamma_2, \gamma_3, \tau_1, \tau_2) = s_1\gamma_1(1 + s_2\gamma_3) - \frac{\gamma_2 + s_2\gamma_3}{s_1\gamma_1}.$$

Here $s_1 = (\tau_1^2 + \tau_2^2)^{-1/2}$ and $s_2 = \tau_1^2\tau_2^2 s_1^4$ are functions of the scales, and the $\hat{\gamma}_i$'s are specified as functions of $y$ under (1.3) and (1.4);

$$(5.3) \qquad \hat{\gamma}_1 \doteq \hat{v} - 2\hat{a}\hat{v}^2, \qquad \hat{\gamma}_2 \doteq \hat{v}(\hat{a} - \hat{d}_1), \qquad \hat{\gamma}_3 \doteq \hat{v}\hat{a}.$$

These $\hat{\gamma}_i$'s are also used to express

$$(5.4) \qquad \hat{z}_2(y) = \hat{\gamma}_1(1 + \hat{\gamma}_3) + \frac{\hat{\gamma}_2}{\hat{\gamma}_1},$$

which is equivalent to (4.9) up to $O(n^{-1/2})$ terms. We calculate $\tilde{\alpha}_2(y, \tau_1, \tau_2)$ for several values of $(\tau_1, \tau_2)$ by the two-step multiscale bootstrap resampling, and fitting the observed $\tilde{z}_2(y, \tau_1, \tau_2) = -\Phi^{-1}(\tilde{\alpha}_2(y, \tau_1, \tau_2))$ to the nonlinear regression model (5.1). Then the estimated $\hat{\gamma}_i$'s are used to compute $\hat{\alpha}_2(y) = \Phi(-\hat{z}_2(y))$ from (5.4).

This procedure is generalized for the three-step multiscale bootstrap resampling. A third-step replicate $Y^{***}$ is generated for each $y^{**}$ by

$$Y^{***} \sim f(y^{***}; y^{**}, \tau_3)$$

using the scale $\tau_3$, and the three-step bootstrap probability is defined by

$$\tilde{\alpha}_3(y, \tau_1, \tau_2, \tau_3) = \Pr\{Y^{***} \in \mathcal{R}; y, \tau_1, \tau_2, \tau_3\}$$
$$= \int \tilde{\alpha}_2(y^*, \tau_2, \tau_3) f(y^*; y, \tau_1) \, dy^*.$$

Then, observed $\tilde{z}_3(y, \tau_1, \tau_2, \tau_3) = -\Phi^{-1}(\tilde{\alpha}_3(y, \tau_1, \tau_2, \tau_3))$ for several values of $(\tau_1, \tau_2, \tau_3)$ are fitted to the nonlinear regression model $\zeta_3$, defined by

$$(5.5) \quad \begin{aligned} &\zeta_3(\gamma_1, \gamma_2, \gamma_3, \gamma_4, \gamma_5, \gamma_6, \tau_1, \tau_2, \tau_3) \\ &= \gamma_1 s_1 (1 + \gamma_3 s_2 + 4\gamma_3^2 s_2^2 + \gamma_5 s_3 + \gamma_6 s_4) \\ &\quad - (\gamma_1 s_1)^{-1}(\gamma_2 + \gamma_3 s_2 + 7\gamma_3^2 s_2^2 + \gamma_4 s_2 + 3\gamma_5 s_3 + 3\gamma_6 s_4), \end{aligned}$$

where $s_1, \ldots, s_4$ are given by

$$s_1 = (\tau_1^2 + \tau_2^2 + \tau_3^2)^{-1/2}, \qquad\qquad s_2 = (\tau_1^2 \tau_2^2 + \tau_2^2 \tau_3^2 + \tau_3^2 \tau_1^2) s_1^4,$$
$$s_3 = (\tau_1^2 \tau_2^2 \tau_3^2 + \tau_2^4 \tau_3^2 + \tau_1^4 (\tau_2^2 + \tau_3^2)) s_1^6, \qquad s_4 = (\tau_1^2 \tau_2^2 \tau_3^2) s_1^6.$$

The least squares estimates for the six $\gamma_i$'s are denoted by $\hat{\gamma}_1, \ldots, \hat{\gamma}_6$. We then compute $\hat{\alpha}_3(y) = \Phi(-\hat{z}_3(y))$ by using the estimated $\hat{\gamma}_i$'s in

$$(5.6) \quad \hat{z}_3(y) = \hat{\gamma}_1 (1 + \hat{\gamma}_3 + 4\hat{\gamma}_3^2 + \hat{\gamma}_6) + \hat{\gamma}_1^{-1}(\hat{\gamma}_2 + \hat{\gamma}_3^2/2 + \hat{\gamma}_4 + \hat{\gamma}_5).$$

Section 7 is mostly devoted to proving the third-order accuracy of $\hat{\alpha}_3(y)$. The justification for the second-order accuracy of $\hat{\alpha}_2(y)$ then immediately follows by ignoring $O(n^{-1})$ terms. As seen in (5.3), $\hat{\gamma}_1$ is $O(1)$, and $\hat{\gamma}_2$ and $\hat{\gamma}_3$ are $O(n^{-1/2})$. The rest of the three $O(n^{-1})$ geometric quantities are defined in Section 7.8. We do not have to know, however, the expressions of the $\hat{\gamma}_i$'s for computing $\hat{\alpha}_3(y)$, because their values are estimated from the nonlinear regression, and the estimation error is only $O(n^{-3/2})$.

It should be noted that there are other asymptotically equivalent expressions for $\zeta_3$ and $\hat{z}_3$ as functions of coefficients transformed from the six $\hat{\gamma}_i$'s; we have shown the two different expressions for $\zeta_2$ and $\hat{z}_2$ as functions of either $\hat{\gamma}_1, \hat{\gamma}_2, \hat{\gamma}_3$ or $\hat{v}, \hat{d}_1, \hat{a}$. The expressions (5.5) and (5.6) are obtained by seeking simple ones.

**6. Examples.** The two procedures in the previous section are applied to the exponential example with parameter values (4.1). By the two-step multiscale bootstrap, the least squares estimates of $\hat{\gamma}_i$'s are

$$\hat{\gamma}_1 = 1.328, \qquad \hat{\gamma}_2 = 0.144, \qquad \hat{\gamma}_3 = 0.137,$$

and the corrected $p$-value is computed as

$$\hat{\alpha}_2(y) = 1 - \Phi\{1.328(1 + 0.137) + \tfrac{0.144}{1.328}\} = 0.0528,$$

which comes closer to the exact $p$-value $\hat{\alpha}_\infty(y) = 0.05$ than $\hat{\alpha}_1(y) = 0.0753$ computed in (4.3). By the three-step multiscale bootstrap, the least squares estimates of the $\hat{\gamma}_i$'s are

$$\hat{\gamma}_1 = 1.328, \qquad \hat{\gamma}_2 = 0.145, \qquad \hat{\gamma}_3 = 0.127,$$

$$\hat{\gamma}_4 = -0.018, \qquad \hat{\gamma}_5 = -0.0004, \qquad \hat{\gamma}_6 = -0.036,$$

and the corrected $p$-value is

$$\hat{\alpha}_3(y) = 1 - \Phi\bigg\{1.328(1 + 0.127 + 0.065 - 0.036)$$
$$+ \frac{0.145 + 0.008 - 0.018 - 0.0004}{1.328}\bigg\} = 0.0509,$$

which is even better than $\hat{\alpha}_2(y) = 0.0528$.

In Table 2 $p$-values are computed for several parameter settings. The bootstrap probabilities are computed numerically ($B = \infty$), but the standard errors due to the bootstrap resampling are shown for $B = 10{,}000$. The first row corresponds to the normal model with (2.2), and the fourth row corresponds to the exponential model with (4.1). The following two rows for each are obtained by changing $n = 10$ to 100 and 1000. Similarly, the last six rows are obtained by changing $\hat{\alpha}_\infty = 0.05$ to 0.95. We observe that all the $p$-values tend to converge to $\hat{\alpha}_\infty$ as $n$ grows, and the corrected $p$-values are faster for convergence than $\hat{\alpha}_0$.

$\tilde{\alpha}_3(y, \tau_1, \tau_2, \tau_3)$ is computed for all the combinations of $(\tau_1, \tau_2, \tau_3)$ values, as noted in the table; five $(\tau_1, 0, 0)$'s, ten $(\tau_1, \tau_2, 0)$'s and twenty $(\tau_1, \tau_2, \tau_3)$'s. Therefore, the numbers of bootstrap probabilities are 5, 15 and 35, respectively, for $\hat{\alpha}_1(y)$, $\hat{\alpha}_2(y)$ and $\hat{\alpha}_3(y)$. The nonlinear regression models are fitted to these bootstrap probabilities, and the least squares estimates of the geometric quantities are calculated; each residual term is weighted inversely proportional to the estimated variance. For stable estimation, ridge regression is also used; a penalty term $\sum_{i=1}^{6} \omega_i \hat{\gamma}_i^2$ with small $\omega_i$ values is added to the residual sum of squares for minimization.

For the exponential distribution, $\hat{\alpha}_k$ is $k$th order accurate ($k = 1, 2, 3$), and, in fact, $|\hat{\alpha}_k - \hat{\alpha}_\infty|$ becomes smaller as $k$ increases in the table. It turns

TABLE 2
*p-values in percent (standard error) for the examples*[*]

| | | | | | | Ridge regression | |
|---|---|---|---|---|---|---|---|
| $n$ | $\hat{\alpha}_0$ | $\hat{\alpha}_{abc}$ | $\hat{\alpha}_1$ | $\hat{\alpha}_2$ | $\hat{\alpha}_3$ | $\hat{\alpha}_2$ | $\hat{\alpha}_3$ |
| | | | Normal distribution ($\hat{\alpha}_\infty = 5.00$) | | | | |
| 10 | 0.85 | 7.75 | 5.29 (0.61) | 5.85 (1.81) | 7.03 (8.04) | 5.67 (1.03) | 6.04 (1.13) |
| 100 | 2.73 | 5.25 | 5.01 (0.37) | 5.05 (1.16) | 5.08 (2.93) | 5.04 (0.78) | 5.06 (0.97) |
| 1000 | 4.12 | 5.03 | 5.00 (0.32) | 5.00 (1.05) | 5.00 (2.22) | 5.00 (0.72) | 5.00 (0.89) |
| | | | Exponential distribution ($\hat{\alpha}_\infty = 5.00$) | | | | |
| 10 | 11.15 | 5.00 | 7.53 (0.31) | 5.28 (0.77) | 5.09 (0.95) | 5.77 (0.60) | 5.13 (0.68) |
| 100 | 6.73 | 5.00 | 5.90 (0.30) | 5.03 (0.94) | 5.01 (1.50) | 5.25 (0.67) | 5.04 (0.81) |
| 1000 | 5.52 | 5.00 | 5.29 (0.30) | 5.00 (0.98) | 5.00 (1.82) | 5.08 (0.69) | 5.01 (0.80) |
| | | | Normal distribution ($\hat{\alpha}_\infty = 95.00$) | | | | |
| 10 | 67.84 | 92.33 | 95.26 (0.18) | 95.20 (0.41) | 95.02 (0.51) | 95.21 (0.34) | 95.07 (0.37) |
| 100 | 90.65 | 94.74 | 95.02 (0.24) | 95.07 (0.84) | 95.09 (1.28) | 95.06 (0.60) | 95.07 (0.70) |
| 1000 | 93.91 | 94.97 | 95.00 (0.28) | 95.00 (0.95) | 95.00 (1.72) | 95.00 (0.67) | 95.00 (0.81) |
| | | | Exponential distribution ($\hat{\alpha}_\infty = 95.00$) | | | | |
| 10 | 98.78 | 95.00 | 97.99 (0.24) | 94.48 (1.31) | 96.12 (7.39) | 95.60 (0.81) | 96.48 (0.56) |
| 100 | 96.49 | 95.00 | 95.95 (0.28) | 94.97 (1.06) | 95.01 (2.71) | 95.24 (0.72) | 95.14 (0.82) |
| 1000 | 95.50 | 95.00 | 95.30 (0.29) | 95.00 (1.02) | 95.00 (2.19) | 95.08 (0.70) | 95.02 (0.81) |

[*]The bootstrap calculation is replaced by integration numerically, and, hence, the number of bootstrap replicates is regarded as $B = \infty$. The standard errors in parentheses are calculated for the case of $B = 10^4$ by the local linearization of the nonlinear regression [Draper and Smith (1998)]. All the combinations of $\tau_1^2 \in \{\frac{10}{3}, \frac{10}{6}, \frac{10}{10}, \frac{10}{15}, \frac{10}{21}\}$, $\tau_2^2 \in \{\frac{10}{6}, \frac{10}{15}\}$, $\tau_3^2 \in \{\frac{10}{6}, \frac{10}{15}\}$ are used for the scales. The total numbers of bootstrap replicates are $5B$, $15B$ and $35B$, respectively, for $\hat{\alpha}_1$, $\hat{\alpha}_2$ and $\hat{\alpha}_3$. For the ridge regression, the penalty weights are $\omega_1 = \omega_2 = 0$ and $\omega_3 = \cdots = \omega_6 = 0.01$.

out that $|\hat{\alpha}_{abc} - \hat{\alpha}_\infty|$ is almost zero here, because $\hat{\alpha}_{abc}$ happens to be third-order accurate for the one-dimensional exponential distribution, as shown in Section 7.7.

For the normal distribution, $\hat{\alpha}_1$, $\hat{\alpha}_2$ and $\hat{\alpha}_3$ are third-order accurate, because $\hat{\gamma}_3 = \cdots = \hat{\gamma}_6 = 0$ under (1.1), as shown in Section 7.8. This may explain why $|\hat{\alpha}_k - \hat{\alpha}_\infty|$ becomes larger as $k$ increases in some of the rows. These four geometric quantities of zero value are estimated from slight differences of bootstrap probabilities, leading to unstable estimation as seen in the large standard errors. This is alleviated by ridge regression; even the worst case in the table $\hat{\alpha}_3 = 6.04 \pm 1.13$ may be allowed in practice. However, the total number of replicates is 350,000 for $\hat{\alpha}_3$, almost comparable to that of the double bootstrap for achieving the same degree of the standard error.

Although $\hat{\alpha}_1$ is first-order accurate for (1.3), it is reasonably accurate even for the exponential model in the table. The total number of replicates is 50,000, yet the standard error is considerably smaller than that of $\hat{\alpha}_3$.

Similar observation holds for the second-order accurate $\hat{\alpha}_2$. The one-step, as well as two-step, multiscale bootstrap may provide a compromise between the number of replicates and the accuracy in practice.

## 7. Asymptotic analysis of the bootstrap methods.

7.1. *A unified approach.* Our approach to assessing the bootstrap methods is not very elegant but rather elementary and brute-force. We explicitly specify a curved coordinate system along $\partial \mathcal{R}$, which is convenient to work on the bootstrap methods. The density function of $Y$ with respect to the curved coordinates is first defined for $\tau = 1$ in Section 7.2 and extended for $\tau > 0$ in Section 7.3. We define a *modified signed distance* by altering $\hat{v}$ slightly, and its distribution function is given in Section 7.4.

It turns out that the $z$-values of the bootstrap probabilities are special cases of the modified signed distance, and our approach gives an asymptotic analysis of the bootstrap methods in a systematic way. Using the result of Section 7.4, a third-order accurate pivot statistic is defined in Section 7.5, and the distribution functions of the bootstrap $z$-values are shown in Sections 7.6 to 7.8, proving the main results of Section 5.

The proofs of lemmas are given in Shimodaira (2004). We have used the computer software *Mathematica* for straightforward and tedious symbolic calculations; the program file is available from the author upon request.

7.2. *Tube-coordinates.* In our curved coordinate system, a point $\eta$ is specified by two parts, a point on $\partial \mathcal{R}$ and the signed distance from it. This is an instance of the coordinate system used for the Weyl tube formula, and we call it tube-coordinates. Below we will define the coordinate system explicitly, and show the expression of the density function of $Y$ in terms of the tube-coordinates. We take an approach similar to that of Kuriki and Takemura (2000).

The density function of the exponential family of distributions is expressed as

$$(7.1) \qquad \exp(\theta^i y_i - \psi(\theta) - h(y)),$$

where $\theta = (\theta^1, \ldots, \theta^p)$ is the natural parameter vector. We denote (7.1) by $f(y; \eta)$ using the expectation parameter vector $\eta = (\eta_1, \ldots, \eta_p) = E(Y)$, the expected value of $Y$. The change of variables $\theta \leftrightarrow \eta$ is one-to-one, and is given by $\eta_i = \partial \psi / \partial \theta^i$, $\theta^i = \partial \phi / \partial \eta_i$, $i = 1, \ldots, p$, where the potential function $\phi(\eta)$ is defined from the cumulant function $\psi(\theta)$ by $\phi(\eta) = \max_\theta \{\theta^i \eta_i - \psi(\theta)\}$. The metric at $\eta$ is denoted as

$$\phi^{ij}(\eta) = \frac{\partial^2 \phi(\eta)}{\partial \eta_i \, \partial \eta_j},$$

and the derivatives of $\phi$ at $\eta = 0$ are denoted as

$$\phi^i = \frac{\partial \phi(\eta)}{\partial \eta_i}\bigg|_0, \qquad \phi^{ij} = \frac{\partial^2 \phi(\eta)}{\partial \eta_i \, \partial \eta_j}\bigg|_0, \qquad \phi^{ijk} = \frac{\partial^3 \phi(\eta)}{\partial \eta_i \, \partial \eta_j \, \partial \eta_k}\bigg|_0, \qquad \text{and so on.}$$

Since the exponential family is not uniquely expressed up to affine transformation, we assume without loss of generality that $\phi^i = 0$ and $\phi^{ij} = \delta_{ij}$, where $\delta_{ij}$ takes value one when $i = j$, otherwise zero. In other words, $E(Y) = 0$ and $\mathrm{cov}(Y)$, the covariance matrix of $Y$, is $I_p$ at $\theta = 0$. We make our asymptotic argument local in a neighborhood of $\eta = 0$ by assuming the local alternatives.

The smooth surface $\partial \mathcal{R}$ of the region $\mathcal{R}$ is specified locally around $\eta = 0$ by

$$\eta_a(u) = u_a, \qquad a = 1, \ldots, p-1; \qquad \eta_p(u) \approx -d^{ab}u_a u_b - e^{abc}u_a u_b u_c,$$

where $u = (u_1, \ldots, u_{p-1})$ is the $(p-1)$-dimensional parameter vector to specify a point $\eta(u)$ on $\partial \mathcal{R}$. $\mathcal{R}$ is specified locally by $\eta_p \leq \eta_p(u)$. It follows from the argument below equation (2.12) of Efron and Tibshirani ([1998]) that $d^{ab} = O(n^{-1/2})$ and $e^{abc} = O(n^{-1})$, and similarly, $\phi^{ijk} = O(n^{-1/2})$ and $\phi^{ijkl} = O(n^{-1})$.

Let $B_i^a(u) = \partial \eta_i / \partial u_a$, $i = 1, \ldots, p$, be the components of a tangent vector of the surface for $a = 1, \ldots, p-1$. They are given explicitly as

$$B_b^a(u) = \delta_{ab}, \qquad b = 1, \ldots, p-1; \qquad B_p^a(u) \approx -2d^{ab}u_b - 3e^{abc}u_b u_c,$$

and the metric in the tangent space is given by

$$\begin{aligned}
\phi^{ab}(u) &= \phi^{ij}(\eta(u)) B_i^a(u) B_j^b(u) \\
&\approx \delta_{ab} + \phi^{abc}u_c \\
&\quad + \{4d^{ac}d^{bd} - 2d^{ac}\phi^{bdp} - 2d^{bd}\phi^{acp} - d^{cd}\phi^{abp} + \tfrac{1}{2}\phi^{abcd}\}u_c u_d,
\end{aligned}$$

(7.2)

where $\phi^{ij}(\eta(u)) \approx \delta_{ij} + \phi^{ija}u_a + \{-d^{ab}\phi^{ijp} + \tfrac{1}{2}\phi^{abij}\}u_a u_b$. Let $B_i^p(u)$, $i = 1, \ldots, p$, be the components of the unit length normal vector orthogonal to the tangent vectors with respect to the metric such that

$$\phi^{ij}(\eta(u)) B_i^a(u) B_j^p(u) = 0, \qquad a = 1, \ldots, p-1;$$

$$\phi^{ij}(\eta(u)) B_i^p(u) B_j^p(u) = 1.$$

The components are calculated explicitly as $B_a^p(u) \approx (2d^{ab} - \phi^{abp})u_b + \{3e^{abc} + d^{ab}\phi^{cpp} + d^{bc}\phi^{app} - 2d^{bd}\phi^{acd} + \phi^{abd}\phi^{cdp} + \tfrac{1}{2}\phi^{abp}\phi^{cpp} - \tfrac{1}{2}\phi^{abcp}\}u_b u_c$, and $B_p^p(u) \approx 1 - \tfrac{1}{2}\phi^{app}u_a + \{-2d^{ac}d^{bc} + \tfrac{1}{2}d^{ab}\phi^{ppp} + \tfrac{1}{2}\phi^{acp}\phi^{bcp} + \tfrac{3}{8}\phi^{app}\phi^{bpp} - \tfrac{1}{4}\phi^{abpp}\}u_a u_b$.

Let $v$ be a scalar, and $(u, v)$ be a $p$-dimensional vector. We consider reparameterization defined by

$$\eta_i(u, v) = \eta_i(u) + B_i^p(u)v, \qquad i = 1, \ldots, p,$$

(7.3)

and assume $\eta \leftrightarrow (u, v)$ is one-to-one at least locally around $\eta = 0$. $(u, v)$ gives the tube-coordinates of the point $\eta$. The boundary $\partial \mathcal{R}$ is expressed simply by $v = 0$, and the region $\mathcal{R}$ is $v \leq 0$. $(u, v)$ is used for indicating the parameter value $\eta = \eta(u, v)$, or the observation $y = \eta(u, v)$. When there is a possibility of confusion, we may write $y \leftrightarrow (\hat{u}, \hat{v})$ instead of $\eta \leftrightarrow (u, v)$.

Since the normal vector is orthogonal to the surface, $\eta(u) = \eta(u, 0) \in \partial \mathcal{R}$ is the projection of $\eta(u, v)$ onto $\partial \mathcal{R}$; $\hat{u}$ is the maximum likelihood estimate under the restricted model specified by $\partial \mathcal{R}$. $\eta(\hat{u}, 0)$ is denoted by $\hat{\eta}(y)$ in Section 1 as a function of $y$. $\hat{v}$ is the signed distance mentioned for (1.1) in Section 3.

$\hat{v}$ is also related to the signed likelihood ratio $R$ [McCullagh (1984) and Severini (2000)] by $R \approx \hat{v} + \frac{1}{6}\hat{\phi}^{ppp}\hat{v}^2 + \{\frac{1}{24}\hat{\phi}^{pppp} - \frac{1}{72}(\hat{\phi}^{ppp})^2\}\hat{v}^3$, where $\hat{\phi}^{ppp}$ and $\hat{\phi}^{pppp}$ are the third and fourth derivatives to the normal direction evaluated at $\eta(\hat{u}, 0)$, instead of $\eta = 0$. This third derivative is associated with the acceleration constant. For the acceleration constant $\hat{a}$, the formula $\hat{a} = -\frac{1}{6}\hat{\phi}^{ppp}$ is obtained directly from equation (2.9) of DiCiccio and Efron (1992), or by using equation (6.7) of Efron (1987) and $\partial^3 \psi / \partial \theta^i \, \partial \theta^j \, \partial \theta^k = -\phi^{ijk}$. The expression for the density function of $(\hat{U}, \hat{V})$ is obtained from $f(y; \eta)$ by change of variables, as shown in the following lemma.

LEMMA 1. *Let $Y \sim f(y; \eta)$ be the exponential family of distributions with $\eta = E(Y)$. Without loss of generality we may assume that $\mathrm{cov}(Y) = I_p$ at $\eta = 0$ and that the true parameter value is specified by $\eta = (0, \dots, 0, \lambda)$ for some $\lambda$, that is, $\eta_a = 0$, $a = 1, \dots, p-1$, $\eta_p = \lambda$, or, equivalently, $u = 0$, $v = \lambda$ using the tube-coordinates $(u, v) \leftrightarrow \eta$. Let $f(\hat{u}, \hat{v}; \lambda)$ be the joint density function of $(\hat{U}, \hat{V}) \leftrightarrow Y$. Then, ignoring the error of $O(n^{-3/2})$, we obtain*

(7.4)
$$\log f(\hat{u}, \hat{v}; \lambda) \approx g(\hat{v}, \lambda) + g^a(\hat{v}, \lambda)\hat{u}_a + g^{ab}(\hat{v}, \lambda)\hat{u}_a\hat{u}_b$$
$$+ g^{abc}(\hat{v}, \lambda)\hat{u}_a\hat{u}_b\hat{u}_c + g^{abcd}(\hat{v}, \lambda)\hat{u}_a\hat{u}_b\hat{u}_c\hat{u}_d,$$

*where the five functions on the right-hand side are defined by $g(\hat{v}, \lambda) = -\frac{1}{2}p\log(2\pi) - \frac{1}{2}(\hat{v} - \lambda)^2 - \frac{1}{8}\phi^{iijj} + \frac{1}{6}(\phi^{ijk})^2 - \frac{1}{3}\phi^{ppp}\lambda^3 - \frac{1}{8}\phi^{pppp}\lambda^4 + \{2d^{aa} - \frac{1}{2}\phi^{aap} + \frac{1}{2}\phi^{ppp} + \frac{1}{2}\phi^{ppp}\lambda^2 + \frac{1}{8}\phi^{pppp}\lambda^3\}\hat{v} + \{-2(d^{ab})^2 + 2d^{ab}\phi^{abp} - \frac{3}{4}(\phi^{abp})^2 - \frac{1}{2}(\phi^{app})^2 - \frac{1}{4}(\phi^{ppp})^2 + \frac{1}{4}\phi^{pppp} + \frac{1}{4}\phi^{aapp}\}\hat{v}^2 - \frac{1}{6}\phi^{ppp}\hat{v}^3 - \frac{1}{24}\phi^{pppp}\hat{v}^4$, $g^a(\hat{v}, \lambda) = \frac{1}{2}\phi^{abb} + \frac{1}{2}\phi^{app}\lambda^2 + \frac{1}{6}\phi^{appp}\lambda^3 + \{-\frac{1}{2}\phi^{app}\lambda - d^{ab}\phi^{bcc} + 5d^{ab}\phi^{bpp} + \phi^{app}d^{bb} - 2\phi^{abc}d^{bc} + \frac{1}{2}\phi^{abp}\phi^{bcc} - \frac{3}{2}\phi^{abp}\phi^{bpp} + \frac{1}{4}\phi^{app}\phi^{bbp} - \frac{3}{4}\phi^{app}\phi^{ppp} + \frac{1}{2}\phi^{abc}\phi^{bcp} - \frac{1}{2}\phi^{abbp} + \frac{1}{2}\phi^{appp} + 6e^{abb} + d^{ab}\phi^{bpp}\lambda^2 - \frac{1}{2}\phi^{abp}\phi^{bpp}\lambda^2 - \frac{1}{4}\phi^{app}\phi^{ppp}\lambda^2\}\hat{v} + \{-d^{ab}\phi^{bpp} + \frac{1}{2}\phi^{abp}\phi^{bpp} + \frac{1}{4}\phi^{app}\phi^{ppp} - \frac{1}{6}\phi^{appp}\}\hat{v}^3$, $g^{ab}(\hat{v}, \lambda) = -\frac{1}{2}\delta_{ab} - d^{ab}\lambda - \frac{1}{2}d^{ab}\phi^{ccp} + \frac{1}{4}\phi^{abcc} - \frac{1}{4}\phi^{acd}\phi^{bcd} + 2d^{ac}d^{bc} - 2d^{ac}\phi^{bcp} - \frac{1}{2}d^{ab}\phi^{ppp}\lambda^2 + \{-d^{ab} + \frac{1}{2}\phi^{abp} - (2d^{ac}d^{bc} - \frac{1}{2}d^{ab}\phi^{ppp} + \frac{1}{4}\phi^{abpp} - \frac{1}{2}\phi^{acp}\phi^{bcp} - \frac{3}{8}\phi^{app}\phi^{bpp})\lambda\}\hat{v}$, $g^{abc}(\hat{v}, \lambda) = -\frac{1}{6}\phi^{abc} - e^{abc}\lambda + \{-2e^{abc} + \frac{1}{3}\phi^{abcp} - \frac{3}{2}d^{ab}\phi^{cpp} + d^{ad}\phi^{bcd} - \frac{1}{2}\phi^{abd}\phi^{cdp} - \frac{1}{4}\phi^{abp}\phi^{cpp}\}\hat{v}$, $g^{abcd}(\hat{v}, \lambda) = -\frac{1}{2}d^{ab}d^{cd} + \frac{1}{2}\phi^{abp}d^{cd} - \frac{1}{24}\phi^{abcd}$.*

7.3. *Changing the scale.* We define a density function $f(y; \eta, \tau)$ with mean $\eta$ and scale $\tau > 0$ by modifying $f(y; \eta)$. Here $\tau$ is regarded as a known constant, whereas $\eta$ is a unknown parameter vector. Let $\phi(\eta, \tau)$ be the potential function of $f(y; \eta, \tau)$, and $\phi(\eta)$ be that for $f(y; \eta)$. Since the density function is defined by specifying the potential function, the following equation gives a definition of $f(y; \eta, \tau)$:

$$(7.5) \qquad\qquad \phi(\eta, \tau) = \phi(\eta)/\tau^2.$$

This $f(y; \eta, \tau)$ comes naturally from the multiscale bootstrap resampling. In fact, the potential function of the replicate $Y^*$ is $\phi(\eta, \tau) = \|\eta\|^2/(2\tau^2)$ for the normal example (2.1) of Section 2, and that is $\phi(\eta, \tau) = -n(1 + \log \eta)/\tau^2$ for the exponential example of Section 4, and thus both agree with (7.5). The same applies to the exponential family, in general, as shown below.

LEMMA 2. *Let $X$ be a $p$-dimensional random vector of the exponential family. We assume that $Y$ is expressed as a sum of $m$ independent $X$'s such that $Y = \sqrt{n}(X_1 + \cdots + X_m)/m$ for $m > 0$, and that the density function is $f(y; \eta)$ when $m = n$. Then $Y \sim f(y; \eta, \tau)$ with $\tau = \sqrt{n/m}$ for $\tau > 0$.*

We continue to use the tube-coordinates defined by the reparameterization $\eta \leftrightarrow (u, v)$ of (7.3). By altering the potential $\phi(\eta, 1)$ to $\phi(\eta, \tau)$, the metric, as well as the tube-coordinates, should have changed if we go back to the specification of $\eta(u)$ and $B^p(u)$ given in the previous section. However, we continue to use the specification with $\tau = 1$ for any $\tau > 0$, so that the reparameterization $\eta \leftrightarrow (u, v)$ does not depend on $\tau$.

LEMMA 3. *Let $f(\hat{u}, \hat{v}; \lambda)$ be the joint density function of $(\widehat{U}, \widehat{V}) \leftrightarrow Y$ given in Lemma 1, and $f(\hat{u}, \hat{v}; \lambda, \tau)$ be that corresponding to $f(y; \eta, \tau)$ with scale $\tau > 0$. Then the expression of $\log f(\hat{u}, \hat{v}; \lambda, \tau)$ is obtained from (7.4) by changing $(\hat{u}, \hat{v})$ to*

$$(7.6) \qquad\qquad \tilde{u} = \hat{u}/\tau, \qquad \tilde{v} = \hat{v}/\tau,$$

*by adding the logarithm of the Jacobian $\log(1/\tau^p)$ to (7.4), and replacing $\phi^{ijk}$, $\phi^{ijkl}$, $d^{ab}$, $e^{abc}$ and $\lambda$, respectively, with*

$$(7.7) \qquad \begin{aligned} \tilde{\phi}^{ijk} &= \tau \phi^{ijk}, & \tilde{\phi}^{ijkl} &= \tau^2 \phi^{ijkl}, \\ \tilde{d}^{ab} &= \tau d^{ab}, & \tilde{e}^{abc} &= \tau^2 e^{abc}, & \tilde{\lambda} &= \lambda/\tau. \end{aligned}$$

7.4. *Modified signed distance.* We consider yet another transformation of the coordinates for expressing the bootstrap $z$-values in modified $\hat{v}$ values. Let $w$ be a scalar variable defined formally by the series

$$(7.8) \qquad\qquad w = v + \sum_{r=0}^{\infty} \bar{c}_r v^r + u_c \sum_{r=0}^{\infty} \bar{b}_r^c v^r,$$

where $v^r$ denotes the $r$th power. The coefficients are $\bar{c}_r = O(n^{-1/2})$ and $\bar{b}_r^c = O(n^{-1})$, and their expressions are specified later. We assume the transformation $(u, v) \leftrightarrow (u, w)$ is one-to-one at least locally around $(u, v) = 0$. By inverting the series in (7.8), we also have

$$(7.9) \qquad v = w - \sum_{r=0}^{\infty} c_r w^r - u_c \sum_{r=0}^{\infty} b_r^c w^r,$$

where $c_r = \bar{c}_r - \sum_{s=0}^{r}(r-s+1)\bar{c}_{r-s+1}\bar{c}_s$, and $b_r^c = \bar{b}_r^c$. The coefficients are $c_r = O(n^{-1/2})$ and $b_r^c = O(n^{-1})$. Let $\widehat{W}$ be the random variable corresponding to $w$; the observed value $\hat{w}$ is defined by (7.8) but using the observed $(\hat{u}, \hat{v})$ instead of $(u, v)$.

We call $\hat{w}$ a modified signed distance characterized by the coefficients $b_r^c$, $c_r$; $\hat{w}$ reduces to $\hat{v}$ when all these coefficients are zero. The $z$-values of the bootstrap probabilities are represented as $\hat{w}$ by appropriately specifying the coefficients. The following lemma plays a key role in studying the distributional properties of the bootstrap probabilities.

LEMMA 4. *Let us assume that the distribution of $Y$ in the tube-coordinates is specified by $(\widehat{U}, \widehat{V}) \sim f(\hat{u}, \hat{v}; \lambda, \tau)$, and the coefficients in (7.9) are of order $b_r^c = O(n^{-1})$ for $r \geq 0$, $c_0 = O(n^{-1/2})$, $c_1 = O(n^{-1})$, $c_2 = O(n^{-1/2})$, $c_3 = O(n^{-1})$ and $c_r = O(n^{-3/2})$ for $r \geq 4$. We define $z_c(\hat{w}; \lambda, \tau)$ from the distribution function of the modified signed distance $\widehat{W}$ as*

$$\Pr\{\widehat{W} \leq \hat{w}\} = \Phi(z_c(\hat{w}; \lambda, \tau)).$$

*Then the $z_c$-formula is, ignoring the error of $O(n^{-3/2})$, expressed as*

$$(7.10) \qquad z_c(\hat{w}; \lambda, \tau) \approx \tau^{-1} g_-(\hat{w}, \lambda) + \tau g_+(\hat{w}, \lambda),$$

*where $g_-(\hat{w}, \lambda) = (\hat{w} - \lambda) - c_0 - \frac{1}{3}\phi^{ppp}\lambda^2 + \frac{1}{6}\phi^{ppp}\lambda\hat{w} + (\frac{1}{6}\phi^{ppp} - c_2)\hat{w}^2 - \frac{1}{6}c_0\phi^{ppp}\lambda - \{c_1 + \frac{1}{3}c_0\phi^{ppp}\}\hat{w} + \{\frac{1}{8}(\phi^{app})^2 + \frac{1}{18}(\phi^{ppp})^2 - \frac{1}{8}\phi^{pppp}\}\lambda^3 + \{-\frac{1}{8}(\phi^{app})^2 + \frac{1}{24}\phi^{pppp}\}\lambda^2\hat{w} + \{-\frac{1}{24}(\phi^{ppp})^2 + \frac{1}{24}\phi^{pppp} - \frac{1}{6}c_2\phi^{ppp}\}\lambda\hat{w}^2 + \{-\frac{1}{72}(\phi^{ppp})^2 + \frac{1}{24}\phi^{pppp} - \frac{1}{3}c_2\phi^{ppp} - c_3\}\hat{w}^3$, and $g_+(\hat{w}, \lambda) = -(d^{aa} + \frac{1}{6}\phi^{ppp}) + \{(d^{ab})^2 - d^{ab}\phi^{abp} + \frac{1}{6}d^{aa}\phi^{ppp} + \frac{1}{2}(\phi^{abp})^2 + \frac{1}{2}(\phi^{app})^2 + \frac{13}{72}(\phi^{ppp})^2 - \frac{1}{4}\phi^{aapp} - \frac{1}{8}\phi^{pppp}\}\hat{w} + \{(d^{ab})^2 - \frac{1}{6}d^{aa}\phi^{ppp} + \frac{1}{8}(\phi^{app})^2 + \frac{5}{72}(\phi^{ppp})^2 - \frac{1}{24}\phi^{pppp}\}\lambda$. Note that the $z_c$-formula does not involve the coefficients $b_r^c$, and that the distribution function of $\widehat{W}$ is characterized by the coefficients $c_r$ with third-order accuracy. The index $c$ of $z_c$ indicates the coefficients $c_r$.*

The true parameter value is assumed to be $(0, \lambda)$ in the $(u, v)$-coordinates for (7.4) and (7.10). If we alter the true parameter value to arbitrary $(u, v)$ with $u \neq 0$, the expression changes as well, and $\Phi^{-1}(\Pr\{\widehat{W} \leq \hat{w}\})$ is denoted as $z_c(\hat{w}; u, v, \tau)$, which reduces to $z_c(\hat{w}; 0, \lambda, \tau) = z_c(\hat{w}; \lambda, \tau)$ when $u = 0$ and $v = \lambda$.

$z_c(\hat{w}; u, v, \tau)$ is used for representing the bootstrap probabilities in particular. The simple bootstrap probability is, for example, $\hat{\alpha}_0(y) = \Pr\{\widehat{V}^* \leq 0; y\} = \Phi(z_c(0; \hat{u}, \hat{v}, 1))$ with all $c_r = 0$. The expression of $z_c(\hat{w}^*; \hat{u}, \hat{v}, \tau)$ is obtained from (7.10) by changing the origin to $\eta(\hat{u})$.

LEMMA 5. *Let $Y^*$ be a replicate of $Y$ distributed conditionally as $Y^* \sim f(y^*; y, \tau)$ with mean $y$ and scale $\tau$, and $\widehat{W}^*$ be the corresponding modified signed distance. Let us denote the conditional distribution of $\widehat{W}^*$ given $y$ as $\Pr\{\widehat{W}^* \leq \hat{w}^*; y\} = \Phi(z_c(\hat{w}^*; \hat{u}, \hat{v}, \tau))$. Then the expression of $z_c(\hat{w}^*; \hat{u}, \hat{v}, \tau)$ is obtained from (7.10) by replacing $\hat{w}$, $\lambda$, $\phi^{ppp}$ and $d_1 = d^{aa}$, respectively, with $\hat{w}^*$, $\hat{v}$,*

$$(7.11) \qquad \hat{\phi}^{ppp} \approx \phi^{ppp} + \{3\phi^{bpp}(2d^{bc} - \phi^{bcp}) - \tfrac{3}{2}\phi^{cpp}\phi^{ppp} + \phi^{cppp}\}\hat{u}_c \quad and$$

$$(7.12) \qquad \hat{d}_1 \approx d^{aa} + \{\tfrac{1}{2}d^{aa}\phi^{cpp} - d^{ab}\phi^{abc} + 3e^{aac}\}\hat{u}_c.$$

*Note that $O(n^{-1})$ terms change only $O(n^{-3/2})$. For example, $d_2 = (d^{ab})^2$ would be replaced with $\hat{d}_2$, but $\hat{d}_2 \approx d_2$.*

### 7.5. *Pivot statistic.*

Although the exactly unbiased $p$-value may not exist in general, a third-order accurate $p$-value can be derived under (1.3) and (1.4). Let $Y^* \sim f(y^*; \hat{\eta}(y), 1)$ be a replicate generated with mean $\hat{\eta}(y)$ instead of $y$, and $\hat{\alpha}_\infty(y)$ be defined as the probability of the corresponding signed distance $\widehat{V}^*$ being greater than or equal to the observed value $\hat{v}$;

$$\hat{\alpha}_\infty(y) = \Pr\{\widehat{V}^* \geq \hat{v}; \hat{\eta}(y)\}.$$

This is the exact $p$-value for the normal example of Section 2 and for the exponential example of Section 4. We will show that $\hat{\alpha}_\infty(y)$ is, in fact, third-order accurate under (1.3) and (1.4).

First, $\hat{z}_\infty(y) = -\Phi^{-1}(\hat{\alpha}_\infty(y))$ is expressed by the $z_c$-formula of Lemma 5. From the definition, $\hat{z}_\infty(y) = z_c(\hat{v}; \hat{u}, 0, 1)$ with all $c_r = 0$ and, thus,

$$(7.13) \quad \begin{aligned} \hat{z}_\infty(y) &\approx \hat{v} - (\hat{d}_1 + \tfrac{1}{6}\hat{\phi}^{ppp}) + \tfrac{1}{6}\hat{\phi}^{ppp}\hat{v}^2 \\ &\quad + \{(d^{ab})^2 - d^{ab}\phi^{abp} + \tfrac{1}{6}d^{aa}\phi^{ppp} \\ &\quad\quad + \tfrac{1}{2}(\phi^{abp})^2 + \tfrac{1}{2}(\phi^{app})^2 + \tfrac{13}{72}(\phi^{ppp})^2 - \tfrac{1}{4}\phi^{aapp} - \tfrac{1}{8}\phi^{pppp}\}\hat{v} \\ &\quad + \{-\tfrac{1}{72}(\phi^{ppp})^2 + \tfrac{1}{24}\phi^{pppp}\}\hat{v}^3. \end{aligned}$$

By comparing (7.13) with (7.8), we find that $\hat{z}_\infty(y)$ can be expressed as $\hat{w}$ with coefficients $\bar{c}_0 = -d^{aa} - \tfrac{1}{6}\phi^{ppp}$, $\bar{c}_1 = (d^{ab})^2 - d^{ab}\phi^{abp} + \tfrac{1}{6}d^{aa}\phi^{ppp} + \tfrac{1}{2}(\phi^{abp})^2 + \tfrac{1}{2}(\phi^{app})^2 + \tfrac{13}{72}(\phi^{ppp})^2 - \tfrac{1}{4}\phi^{aapp} - \tfrac{1}{8}\phi^{pppp}$, $\bar{c}_2 = \tfrac{1}{6}\phi^{ppp}$, $\bar{c}_3 = -\tfrac{1}{72}(\phi^{ppp})^2 +$

$\frac{1}{24}\phi^{pppp}$, $\bar{b}_0^c = -\frac{1}{2}d^{aa}\phi^{cpp} + d^{ab}\phi^{abc} - 3e^{aac}$ and $\bar{b}_2^c = \frac{1}{2}\phi^{bpp}(2d^{bc} - \phi^{bcp}) - \frac{1}{4}\phi^{cpp}\phi^{ppp} + \frac{1}{6}\phi^{cppp}$. Then the distribution function of $\hat{z}_\infty(y)$ is obtained immediately from Lemma 4 as shown below.

LEMMA 6. *Let us consider a statistic*

$$\hat{z}_q(y) \approx \hat{z}_\infty(y) + q_0 + q_1\hat{v} + q_2\hat{v}^2 + q_3\hat{v}^3 + \hat{u}_c g^c(\hat{v}),$$

*where the coefficients are* $q_0 = O(n^{-1/2})$, $q_1 = O(n^{-1})$, $q_2 = O(n^{-1/2})$ *and* $q_3 = O(n^{-1})$, *and* $g^c(\hat{v}) = O(n^{-1})$, $c = 1, \ldots, p - 1$, *representing arbitrary polynomials of* $\hat{v}$. *The index* $q$ *of* $z_q$ *indicates the coefficients. Assuming* $(\hat{U}, \hat{V}) \sim f(\hat{u}, \hat{v}; \lambda, 1)$, *the distribution function of* $\hat{z}_q(y)$ *is expressed as*

$$\Pr\{\hat{z}_q(Y) \le x; \lambda\}$$

$$
\begin{aligned}
(7.14) \quad &\approx \Phi[x - \lambda - q_0 - \tfrac{1}{3}\phi^{ppp}\lambda^2 + \tfrac{1}{6}\phi^{ppp}\lambda x - q_2 x^2 \\
&\quad + \{(d^{ab})^2 + \tfrac{1}{8}(\phi^{app})^2 + \tfrac{7}{72}(\phi^{ppp})^2 - \tfrac{1}{24}\phi^{pppp} - \tfrac{1}{6}\phi^{ppp}q_0\}\lambda \\
&\quad + \{-q_1 - 2q_2(d^{aa} + \tfrac{1}{6}\phi^{ppp} - q_0)\}x + \{-\tfrac{1}{8}(\phi^{app})^2 + \tfrac{1}{24}\phi^{pppp}\}\lambda^2 x \\
&\quad + \{\tfrac{1}{3}\phi^{ppp}q_2 + 2q_2^2 - q_3\}x^3 + \{\tfrac{1}{8}(\phi^{app})^2 + \tfrac{1}{18}(\phi^{ppp})^2 - \tfrac{1}{8}\phi^{pppp}\}\lambda^3 \\
&\quad\quad\quad\quad\quad + \{-\tfrac{5}{72}(\phi^{ppp})^2 + \tfrac{1}{24}\phi^{pppp} - \tfrac{1}{6}\phi^{ppp}q_2\}\lambda x^2].
\end{aligned}
$$

For $\lambda = 0$, the distribution function is $\Pr\{\hat{z}_q(Y) \le x; 0\} \approx \Phi[x - q_0 - q_2 x^2 + \{-q_1 - 2q_2(d^{aa} + \tfrac{1}{6}\phi^{ppp} - q_0)\}x + \{\tfrac{1}{3}\phi^{ppp}q_2 + 2q_2^2 - q_3\}x^3]$. In particular, $\Pr\{\hat{z}_\infty(Y) \le x; 0\} \approx \Phi(x)$ and, thus, $\hat{z}_\infty(y)$ is a third-order accurate pivot statistic. We obtain $\Pr\{\hat{\alpha}_\infty(Y) < \alpha; \eta\} \approx \alpha$ for $\eta \in \partial\mathcal{R}$, proving the third-order accuracy of $\hat{\alpha}_\infty(y)$.

The reverse of the above statement also holds. $\hat{\alpha}_q(y) = \Phi(-\hat{z}_q(y))$ is a third-order accurate $p$-value if and only if $q_0 \approx q_1 \approx q_2 \approx q_3 \approx 0$. If we confine our attention to $\hat{\alpha}_q(y)$ defined only from $\hat{v}$ and the geometric quantities $d^{ab}$, $e^{aac}$, $\phi^{ij}$, $\phi^{ijk}$ and $\phi^{ijkl}$ evaluated at $\hat{\eta}(y)$, then $\hat{u}_c g^c(\hat{v})$ in $\hat{z}_q(y)$ comes only from $q_r$'s by the replacements shown in Lemma 5. Thus, $\hat{\alpha}_q(y)$ is a third-order accurate $p$-value if and only if $\hat{\alpha}_q(y) \approx \hat{\alpha}_\infty(y)$. Similarly, $\hat{\alpha}_q(y)$ is second-order accurate if and only if $q_0 \doteq q_2 \doteq 0$ and, thus, $\hat{\alpha}_q(y) \doteq \hat{\alpha}_\infty(y)$.

$\hat{z}_\infty(y)$ is equivalent to other pivots in the literature up to $O(n^{-1})$ terms. Under (1.1) and (1.4), $\phi^{ijk} = \phi^{ijkl} = 0$ and, thus, (7.13) reduces to $\hat{z}_\infty(y) \approx \hat{v} - \hat{d}_1 + \hat{d}_2\hat{v}$, giving (3.8), the pivot of Efron (1985). Under (1.3), the modified signed likelihood ratio [Barndorff-Nielsen (1986) and Barndorff-Nielsen and Cox (1994)] has been known as a third-order accurate pivot, and it is expressed as $R^* = R + (1/R)\log(U/R)$ in the notation of Severini [(2000), page 251], where $U$ is defined using the log-likelihood derivatives. A straightforward calculation shows that $U \approx \hat{v} - \hat{d}_1\hat{v}^2 + \{\tfrac{1}{2}(d^{aa})^2 + d^{ab}d^{ab} - \tfrac{1}{4}\phi^{aapp} - d^{ab}\phi^{abp} + \tfrac{1}{2}(\phi^{abp})^2 + \tfrac{1}{2}(\phi^{app})^2 + \tfrac{1}{8}(\phi^{ppp})^2 - \tfrac{1}{12}\phi^{pppp}\}\hat{v}^3$, and that $R^* \approx \hat{z}_\infty(y)$ in the moderate deviation region.

7.6. *Accuracy of the bootstrap probability.* Since the event $Y^* \in \mathcal{R}$ is equivalent to the event $\hat{V}^* \leq 0$, the $z$-value of the bootstrap probability with scale $\tau$ is expressed by the $z_c$-formula of Lemma 5; $\tilde{z}_1(y, \tau) = -z_c(0; \hat{u}, \hat{v}, \tau)$ with all $c_r = 0$. From (7.10), we obtain a refined version of (4.8), erring only $O(n^{-3/2})$,

$$
\tilde{z}_1(y, \tau) \approx \tau^{-1}[\hat{v} + \tfrac{1}{3}\hat{\phi}^{ppp}\hat{v}^2 - \{\tfrac{1}{8}(\phi^{app})^2 + \tfrac{1}{18}(\phi^{ppp})^2 - \tfrac{1}{8}\phi^{pppp}\}\hat{v}^3]
$$

$$
(7.15) \qquad + \tau[(\hat{d}_1 + \tfrac{1}{6}\hat{\phi}^{ppp})
$$

$$
- \{(d^{ab})^2 - \tfrac{1}{6}d^{aa}\phi^{ppp} + \tfrac{1}{8}(\phi^{app})^2 + \tfrac{5}{72}(\phi^{ppp})^2 - \tfrac{1}{24}\phi^{pppp}\}\hat{v}].
$$

It follows from (7.15) that $\tau\tilde{z}_1(y, \tau)$ is expressed as $\hat{w}$ and, thus, $\tau\tilde{z}_1(y, \tau) \approx \hat{z}_q(y)$ by choosing the coefficients appropriately. They are $c_0 = (d^{aa} + \tfrac{1}{6}\phi^{ppp})\tau^2$, $c_1 = (-(d^{ab})^2 - \tfrac{1}{2}d^{aa}\phi^{ppp} - \tfrac{1}{8}(\phi^{app})^2 - \tfrac{13}{72}(\phi^{ppp})^2 + \tfrac{1}{24}\phi^{pppp})\tau^2$, $c_2 = \tfrac{1}{3}\phi^{ppp}$, and $c_3 = -\tfrac{1}{8}(\phi^{app})^2 - \tfrac{5}{18}(\phi^{ppp})^2 + \tfrac{1}{8}\phi^{pppp}$ for $\hat{w}$, or, equivalently, $q_0 = (1 + \tau^2)(d^{aa} + \tfrac{1}{6}\phi^{ppp})$, $q_1 = -(1 + \tau^2)(d^{ab})^2 + d^{ab}\phi^{abp} + \tfrac{1}{4}\phi^{aapp} - \tfrac{1}{2}(\phi^{abp})^2 - \tfrac{1}{8}(4 + \tau^2)(\phi^{app})^2 + \tfrac{1}{6}(-1 + \tau^2)d^{aa}\phi^{ppp} - \tfrac{1}{72}(13 + 5\tau^2)(\phi^{ppp})^2 + \tfrac{1}{24}(3 + \tau^2)\phi^{pppp}$, $q_2 = \tfrac{1}{6}\phi^{ppp}$, $q_3 = -\tfrac{1}{8}(\phi^{app})^2 - \tfrac{1}{24}(\phi^{ppp})^2 + \tfrac{1}{12}\phi^{pppp}$ for $\hat{z}_q(y)$. The distribution function of $\tau\tilde{z}(y, \tau)$ is obtained from (7.10) or (7.14). In particular, the distribution function of $\hat{z}_0(y) = \tilde{z}_1(y, 1)$ under $\lambda = 0$, $\tau = 1$ is

$$
\Pr\{\hat{z}_0(Y) \leq x; 0\}
$$

$$
\approx \Phi[x - (2d^{aa} + \tfrac{1}{3}\phi^{ppp}) - \tfrac{1}{6}\phi^{ppp}x^2
$$

$$
(7.16) \qquad + \{2(d^{ab})^2 - d^{ab}\phi^{abp} + \tfrac{1}{3}d^{aa}\phi^{ppp} + \tfrac{1}{2}(\phi^{abp})^2
$$

$$
+ \tfrac{5}{8}(\phi^{app})^2 + \tfrac{11}{36}(\phi^{ppp})^2 - \tfrac{1}{4}\phi^{aapp} - \tfrac{1}{6}\phi^{pppp}\}x
$$

$$
+ \{\tfrac{11}{72}(\phi^{ppp})^2 + \tfrac{1}{8}(\phi^{app})^2 - \tfrac{1}{12}\phi^{pppp}\}x^3],
$$

showing the first-order accuracy of $\hat{\alpha}_0(y)$.

Remark A of Efron and Tibshirani (1998) discusses a calibrated bootstrap probability, denoted $\hat{\alpha}_{\text{double}}(y)$ here, using the double bootstrap of Hall (1992). Similarly to the two-level bootstrap, thousands of $Y^*$ are generated around $\hat{\eta}(y)$. Then $\hat{\alpha}_0(y^*)$ is computed for each $y^*$. The expression of $\hat{z}_{\text{double}}(y) = \Phi^{-1}[\Pr\{\hat{z}_0(Y^*) \leq \hat{z}_0(y); \hat{\eta}(y)\}]$ is obtained from (7.16) by the replacements of Lemma 5, and a straightforward calculation shows that $\hat{z}_{\text{double}}(y) \approx \hat{z}_\infty(y)$, proving the third-order accuracy of $\hat{\alpha}_{\text{double}}(y)$.

7.7. *Accuracy of the two-level bootstrap.* The expression of $\hat{z}_0(y)$ is obtained from (7.15) by letting $\tau = 1$, and $\hat{z}_0(\hat{\eta}(y)) \approx \hat{d}_1 + \tfrac{1}{6}\hat{\phi}^{ppp}$ is obtained from it by letting $\hat{v} = 0$. By substituting these expressions, as well as $\hat{a} = -\tfrac{1}{6}\hat{\phi}^{ppp}$ for those in (2.3), we find that $\hat{z}_{\text{abc}}(y)$ is expressed as $\hat{w}$, or, equivalently, $\hat{z}_q(y)$ with coefficients $q_0 = q_2 = 0$, $q_1 = -2(d^{ab})^2 + \tfrac{1}{4}\phi^{aapp} + d^{ab}\phi^{abp} -$

$\frac{1}{2}(\phi^{abp})^2 - \frac{5}{8}(\phi^{app})^2 - \frac{1}{4}(\phi^{ppp})^2 + \frac{1}{6}\phi^{pppp}$ and $q_3 = -\frac{1}{8}(\phi^{app})^2 - \frac{1}{8}(\phi^{ppp})^2 + \frac{1}{12}\phi^{pppp}$. The distribution function is then obtained from Lemma 6. For $\lambda = 0$, it becomes

$$(7.17) \qquad \Pr\{\hat{z}_{\mathrm{abc}}(Y) \leq x; 0\} \approx \Phi(x - q_1 x - q_3 x^3),$$

showing the second-order accuracy of $\hat{\alpha}_{\mathrm{abc}}(y)$.

For the exponential example of Section 4, $p = 1$, $\phi^{111} = -2/\sqrt{n}$, $\phi^{1111} = 6/n$ and all the other quantities in $q_1$ and $q_3$ are zero. Therefore, $q_1 = q_3 = 0$, and $\hat{z}_{\mathrm{abc}}(y)$ turns out to be third-order accurate, explaining the high accuracy of $\hat{\alpha}_{\mathrm{abc}}(y)$ observed in Table 2.

7.8. *Accuracy of the multistep-multiscale bootstrap.* Using the expressions (7.4) and (7.15), the expression of $\tilde{z}_2(y, \tau_1, \tau_2)$ is obtained by the integration

$$(7.18) \qquad \tilde{z}_2(y, \tau_1, \tau_2) = \Phi^{-1}\left\{\int \Phi(\tilde{z}_1(y^*, \tau_2))f(y^*; y, \tau_1)\, dy^*\right\}.$$

By repeating the same integration using $\tilde{z}_2(y^*, \tau_2, \tau_3)$ instead of $\tilde{z}_1(y^*, \tau_2)$, we obtain the expression of $\tilde{z}_3(y, \tau_1, \tau_2, \tau_3)$ as given below.

LEMMA 7. *Let us define the following six geometric quantities using the derivatives evaluated at $\eta = 0$: $\gamma_1 = \lambda + \frac{1}{3}\lambda^2\phi^{ppp} + \lambda^3\{-\frac{1}{8}(\phi^{app})^2 - \frac{1}{18}(\phi^{ppp})^2 + \frac{1}{8}\phi^{pppp}\}$, $\gamma_2 = \lambda\{-d^{aa} - \frac{1}{6}\phi^{ppp}\} + \lambda^2\{(d^{ab})^2 - \frac{1}{2}d^{aa}\phi^{ppp} + \frac{1}{8}(\phi^{app})^2 + \frac{1}{72}(\phi^{ppp})^2 - \frac{1}{24}\phi^{pppp}\}$, $\gamma_3 = -\frac{1}{6}\lambda\phi^{ppp} + \lambda^2\{\frac{1}{4}(\phi^{app})^2 + \frac{1}{9}(\phi^{ppp})^2 - \frac{1}{8}\phi^{pppp}\}$, $\gamma_4 = \lambda^2\{-d^{ab}\phi^{abp} + \frac{1}{3}d^{aa}\phi^{ppp} + \frac{1}{2}(\phi^{abp})^2 + \frac{1}{2}(\phi^{app})^2 + \frac{2}{9}(\phi^{ppp})^2 - \frac{1}{4}\phi^{aapp} - \frac{1}{6}\phi^{pppp}\}$, $\gamma_5 = \lambda^2\{-\frac{1}{8}(\phi^{app})^2 - \frac{1}{8}(\phi^{ppp})^2 + \frac{1}{12}\phi^{pppp}\}$ and $\gamma_6 = \lambda^2\{-\frac{1}{8}(\phi^{app})^2 - \frac{1}{8}(\phi^{ppp})^2 + \frac{1}{24}\phi^{pppp}\}$. Those evaluated at $\hat{\eta}(y)$, denoted $\hat{\gamma}_1, \ldots, \hat{\gamma}_6$, are obtained by replacing $\lambda$, $\phi^{ppp}$ and $d^{aa}$, respectively, with $\hat{v}$, (7.11) and (7.12) as shown in Lemma 5. Then we have*

$$(7.19) \qquad \tilde{z}_3(y, \tau_1, \tau_2, \tau_3) \approx \zeta_3(\hat{\gamma}_1, \hat{\gamma}_2, \hat{\gamma}_3, \hat{\gamma}_4, \hat{\gamma}_5, \hat{\gamma}_6, \tau_1, \tau_2, \tau_3)$$

*using the $\zeta_3$-function of (5.5). Since (7.19) errs only $O(n^{-3/2})$ for any values of $(\tau_1, \tau_2, \tau_3)$, the nonlinear regression for three-step multiscale bootstrap probabilities in Section 5 estimates $\hat{\gamma}_i$'s up to $O(n^{-1})$ terms.*

If we define $\hat{z}_3(y)$ of (5.6) using the $\hat{\gamma}_i$'s defined above, we can easily verify

$$(7.20) \qquad \hat{z}_3(y) \approx \hat{z}_\infty(y)$$

by comparing (5.6) with (7.13). This proves the third-order accuracy of $\hat{\alpha}_3(y)$ under (1.3) and (1.4).

For the multivariate normal model of (1.1), $\phi(\eta) = \|\eta\|^2/2$ and, thus, $\phi^{ijk} = \phi^{ijkl} = 0$. This implies $\gamma_3 = \cdots = \gamma_6 = 0$, proving the third-order accuracy of $\hat{\alpha}_1(y)$ and $\hat{\alpha}_2(y)$ under (1.1) and (1.4).

**Acknowledgments.** I wish to thank the referees and the Associate Editor who handled this article for their very helpful constructive suggestions. The earlier version of the manuscript was prepared during my stay at Stanford University arranged by Brad Efron.

Department of Mathematical
  and Computing Sciences
Tokyo Institute of Technology
2-12-1 Ookayama, Meguro-ku
Tokyo 152-8552
Japan
E-MAIL: shimo@is.titech.ac.jp
URL: www.is.titech.ac.jp/˜shimo/